\documentclass[11pt]{article}
\usepackage{amssymb}

   \csname @addtoreset\endcsname{equation}{section}
    \csname @addtoreset\endcsname{figure}{section}
    \csname @addtoreset\endcsname{table}{section}

\newcounter{thanksnum}
\def\thanksnumber#1
{\setcounter{thanksnum}{\value{footnote}}\setcounter{footnote}{#1}%
                     \addtocounter{footnote}{-1}\footnotemark
                     \setcounter{footnote}{\value{thanksnum}}}
\def\newtheoremz#1{\@ifnextchar[{\@othmz{#1}}{\@nthmz{#1}}}

\def\@nthmz#1#2{%
\@ifnextchar[{\@xnthmz{#1}{#2}}{\@ynthmz{#1}{#2}}}

\def\@xnthmz#1#2[#3]{\expandafter\@ifdefinable\csname #1\endcsname
{\@definecounter{#1}\@addtoreset{#1}{#3}%
\expandafter\xdef\csname the#1\endcsname{\expandafter\noexpand
  \csname the#3\endcsname \@thmcountersepz \@thmcounterz{#1}}%
\global\@namedef{#1}{\@thmz{#1}{#2}}\global\@namedef{end#1}{\@endtheoremz}}}

\def\@ynthmz#1#2{\expandafter\@ifdefinable\csname #1\endcsname
{\@definecounter{#1}%
\expandafter\xdef\csname the#1\endcsname{\@thmcounterz{#1}}%
\global\@namedef{#1}{\@thm{#1}{#2}}\global\@namedef{end#1}{\@endtheoremz}}}

\def\@othmz#1[#2]#3{\expandafter\@ifdefinable\csname #1\endcsname
  {\global\@namedef{the#1}{\@nameuse{the#2}}%
\global\@namedef{#1}{\@thmz{#2}{#3}}%
\global\@namedef{end#1}{\@endtheoremz}}}

\def\@thmz#1#2{\refstepcounter
    {#1}\@ifnextchar[{\@ythmz{#1}{#2}}{\@xthmz{#1}{#2}}}

\def\@xthmz#1#2{\@begintheoremz{#2}{\csname the#1\endcsname}\ignorespaces}
\def\@ythmz#1#2[#3]{\@opargbegintheoremz{#2}{\csname
       the#1\endcsname}{#3}\ignorespaces}

\def\@thmcounterz#1{\noexpand\arabic{#1}}
\def\@thmcountersepz{.}
\def\@begintheoremz#1#2{ \trivlist \item[\hskip \labelsep{\bf #1\ #2}]}
\def\@opargbegintheoremz#1#2#3{ \trivlist
      \item[\hskip \labelsep{\bf #1\ #2\ (#3)}]}
\def\@endtheoremz{\endtrivlist}

\newtheorem{theorem}{Theorem}[section]
\newtheorem{lemma}{Lemma}[section]

\newtheorem{proposition}{Proposition}[section]
\newtheorem{corollary}{Corollary}[section]
\newtheorem{condition}{Condition}[section]
\newtheorem{definition}{Definition}[section]
\newtheorem{remark}{Remark}[section]

\def\e{\varepsilon}

\def\defi{\stackrel{{\scriptscriptstyle \Delta}}{=}}

\def\a{\alpha}
\def\d{\delta}
\def\o{\omega}
\def\O{\Omega}

\def\F{{\cal F}}
\def\w{\widehat}
\def\ww{\widetilde}
\def\Ind{{\mathbb{I}}}

\def\esssup{\mathop{\rm ess\, sup}}
\def\const{{\rm const\,}}

\def\R{{\bf R}}
\def\E{{\bf E}}
\def\P{{\bf P}}

\def\Z{{\cal Z}}

\def\L{L}

\def\b{\beta}
\def\s{\delta}
\def\g{\gamma}

\def\W{{\cal W}^*}
\def\X{{\cal X}}
\def\t{\theta}
\def\oo{\bar}
\def\s{\sigma}

\def\p{\partial}

\def\A{{\cal A}}
\def\M{{\cal M}}

\def\L{{\cal L}}
\def\I{{\, \cal I}}

\newcommand{\be}{\begin{equation}}
\newcommand{\ee}{\end{equation}}
\newcommand{\bd}{\begin{displaymath}}
\newcommand{\ed}{\end{displaymath}}
\newcommand{\ba}{\begin{array}{ll}}
\newcommand{\ea}{\end{array}}
\newcommand{\baa}{\begin{eqnarray}}
\newcommand{\eaa}{\end{eqnarray}}
\newcommand{\baaa}{\begin{eqnarray*}}
\newcommand{\eaaa}{\end{eqnarray*}}   \font\sm=cmr10

\def\W{{\cal W}}

\def\CC{{\cal C}}

\def\LU{L}

\def\I{{ \d}}

\date{First version:  23 Jun 2006. Revised version:  27 July 2010}
\title{Representation of
functionals  of Ito processes and their first exit times
\footnote{Accepted to {\em Stochastics}.}}
\author{
Nikolai Dokuchaev\\ {\sm  Department of Mathematics \& Statistics,
Curtin University of Technology},\\ {\sm GPO Box U1987, Perth, 6845
Western Australia. Email N.Dokuchaev at curtin.edu.au}}
\begin{document}
\maketitle
\begin{abstract} The representation theorem is obtained for functionals of non-Markov
processes and their first exit times from  bounded domains. These
functionals are represented via solutions of backward parabolic Ito
equations. As an example of applications,  analogs of forward
Kolmogorov equations are derived for conditional probability density
  functions of Ito processes being killed on the
  boundary. In addition,
a maximum principle and a contraction property are established  for
SPDEs in bounded domains.
\end{abstract}
\section{Introduction}
In the present paper, we study representation of integrals of
stochastic non-Markov processes and their first exit times via
stochastic partial differential equations. It is a generalization of
the classical Kolmogorov representation for Markov diffusion
processes.
\par
Let  a region $D\subset \R^n$ be given, let $T>0$ be a  terminal
time, let $\F_t$ be a filtration, and let $y^{x,s}(t)$ be an Ito
process adapted to $\F_t$ and such that $y^{x,s}(s)=x$, $x\in D$,
$s<T$. Further, let  $\tau^{x,s}$ be the first exit time from
$D\times[0,T)$ for the vector $(y^{x,s}(t),t)$, and let $\Psi$ and
$\xi$ be some functions. Our goal is to represent conditional
expectations
 \be\label{rep0}
\ww p(x,s,\o)\defi \E\Big\{\Psi(y^{x,s}(T))\Ind_{\{T\le \tau^{x,s}
\}}\,|\,\F_s\Big\}+\E\biggl\{\int_s^{\tau^{x,s}}
\,\xi(y^{x,s}(t),t,\o)\,dt\,\bigr|\,\F_s\biggr\} \ee as the
solutions of boundary value problems for  stochastic partial
differential equations. This representation has many important
applications. In particular, the representation via solution of a
SPDE helps to establish some regularity properties for $\ww p$ and
$\tau^{x,s}$, since there is certain regularity for the solutions of
SPDEs.
\par
For the representation, we will use backward parabolic Ito
  equations, i.e., the equations with Cauchy condition
 at terminal time $t=T$. These equations are  analogs of Kolmogorov backward equations
 for non-Markov processes. We will also consider  forward parabolic Ito
  equations, i.e., the equations with Cauchy condition
 at initial time; they can be regarded as analogs of forward
 Kolmogorov equations.
\par
  Boundary
value problems for forward parabolic Ito equations were intensively
studied; see, e.g., Al\'os et al (1999), Bally {\it et al} (1994),
Chojnowska-Michalik and Goldys (1995), Da Prato and Tubaro (1996),
Gy\"ongy (1998), Kim (2004), Krylov (1999), Maslowski (1995),
Pardoux (1993),
 Rozovskii (1990), Walsh (1986), Zhou (1992), the author's papers
 (1995), (2005),
and the bibliography there.  Note that the difference between
backward and forward equations is not that important for the
deterministic equations because one can always make a change of time
variable and convert a backward equation to a forward one and
opposite. But it cannot be done so easily for stochastic equations,
because the solution needs to be adapted to the driving Brownian
motion.  Therefore, backward stochastic partial differential
equations with boundary conditions at final time require special
consideration.  A possible approach is to consider so-called
Ito-Bismut backward equations when the diffusion term is not given a
priori but has to be found. These backward SPDEs were also widely
studied; see, e.g., Pardoux and Peng
 (1990), Hu and Peng  (1991),  Dokuchaev (1992),
 (2003),(2010),
  Yong and Zhou
 (1999),
 Pardoux and Rascanu (1998), Ma and Yong (1999), Hu {\it et al} (2002), Confortola
 (2007), and references here.  The duality between linear forward and
backward equations was  studied by Zhou (1992) for a domain without
boundary, and by the author (1992) for the domains with boundaries.
A different type of backward equations was described in Chapter 5 of
Rozovskii (1990).
\par
The representation of expectations (\ref{rep0}) via SPDEs was
established before for the following cases:
\begin{itemize}
\item For the classical Markovian setting then $y^{x,s}(t)$ is a diffusion Markov processes;
\item  For the case of non-Markov $y^{x,s}(t)$ in the  entire space, i.e.,  when $D=\R^n$, i.e., for the problem
 without random first exit times.
\end{itemize}
\par
The known representation theorems for non-Markov processes in
$D=\R^n$  was never extended on the case of domains with boundary.
Let us explain why it is non-trivial.
\par
 The main difficulty in the implementation of this approach to
 the non-Markov Ito processes and the related SPDEs is the following.  One  needs again a priori certain
smoothness for the solution $p(\cdot)$
 of a backward SPDE,  to apply Ito-Ventsell
formula for the process $p(y^{x,s}(t,\o),t,\o)$.
 However, the previously known results about
regularity of the solution of the backward SPDE for $p$ were
insufficient for the case of domains with boundary. Therefore, the
representation result was never obtained for this case.
Correspondingly, it was unknown if the forward parabolic Ito
equation for the conditional density of a non-Markov process in the
entire space can be used for the process being killed on the
boundary, given additional
 Dirichlet boundary value condition on this boundary. As far as we know, the first
attempt to solve it was made in the author's paper (1992) for a very
special case. In the present paper, we have proved this fact
together with representation (\ref{rep0}) for some $p$ derived from
a backward parabolic Ito equation (Theorems \ref{Threpr} and Theorem
\ref{Thprobd}).
\par
The present paper uses the additional regularity in the form of the
so-called second fundamental inequality (Theorem \ref{Th2FE}): the
solution $(p,\chi)$
 of the backward equation
 has $L_2$-integrable second derivatives for $p$ and the first derivatives for $\chi$. This
additional regularity of the solutions of the backward equations
appears to be sufficient to obtain the representation theorem.  To
ensure this regularity, we required  additional Condition
\ref{condA} which is  a strengthened version of the standard
coercivity condition (Condition \ref{cond3.1.A}). We emphasize that,
without this new condition, representation theorem for (\ref{rep0})
is still not established, and an equation for the probability
density function of the Ito process being killed on the boundary is
still unknown (even if it easy to believe that one can use the SPDE
for the density from the case of entire domain with additional the
Dirichlet condition imposed on the boundary).
\par
As a corollary, we obtained  the equation for the conditional
probability density function of an Ito process being killed on the
boundary of a domain (Theorem \ref{Thprobd}). This is a new result
even given that the corresponding result for entire domain was known
for a long time (see, e.g., Theorem 5.3.1 from Rozovskii (1990)). As
an additional corollary, we obtained the "maximum principle": the
solution of the forward or backward equation in the cylinder
$D\times[0,T]$ is nonnegative if the free terms are nonnegative.
Further, we proved that the dynamic of the homogeneous equations is
of the contraction type: $\E\int_D|u(x,T,\o)|dx\le
\E\int_D|u(x,0,\o)|dx$ for the solutions of the forward equations,
and $\esssup_{x,\o} |p(x,t,\o)|\le\esssup_{x,\o} |p(x,T,\o)|$ for
the solutions of the backward equations.
 (Theorems \ref{Th6.1}-\ref{Th6.4}).
\par
 The  paper is organized as follows. In Section two we collect
notation and definitions. Sections three contains some facts about
the regularity of SPDEs, including the second fundamental inequality
for backward equations. In Section four, the main result is
presented. The proof of this result is given in Section five.
Sections six and seven contain applications.
\section{Definitions}\label{SecD}
\subsection{Spaces and classes of functions.} 
We a given an open domain $D\subseteq\R^n$ such that either $D=\R^n$
or $D$ is bounded  with $C^{2+\a}$-smooth boundary $\p D$ for some
$\a>0$; if $n=1$, then the condition of smoothness is not required.
Let $T>0$ be given, and let $Q\defi D\times (0,T)$.
\par We are given a standard complete probability space $(\O,\F,\P)$
and a right-continuous filtration $\F_t$ of complete $\s$-algebras
of events, $t\ge 0$; we denote by $\o$ the elements of the set
$\O=\{\o\}$. We are also given  a $N$-dimensional process
$w(t)=(w_1(t),...,w_N(t))$ with independent components such that it
is a Wiener process with respect to $\F_t$.
\par
We denote by $\|\cdot\|_{ X}$ the norm in a linear normed space
$X$, and
 $(\cdot, \cdot )_{ X}$ denotes the scalar product in  a Hilbert space $
X$.
\par
 We denote Euclidean norm in $\R^k$ as $|\cdot|$, and $\bar G$ denotes
the closure of a region $G\subset\R^k$.
\par We introduce some spaces
of real valued functions.
\par
We denote by ${W_q^m}(D)$  the Sobolev  space of functions that
belong to $L_q(D)$ together with first $m$ derivatives, $q\ge 1$. In
particular, $$ \|u\|_{W_2^1(D)}\defi
\left(\|u\|_{L_2(D)}^2+\sum_{i=1}^n\left\|\frac{\p u}{\p
x_i}\right\|_{L_2(D)}^2 \right)^{1/2}.
$$
\par Let $H^0\defi L_2(D)$,
and let $H^1\defi \stackrel{\scriptscriptstyle 0}{W_2^1}(D)$ be the
closure in the ${W}_2^1(D)$-norm of the set of all smooth functions
$u:D\to\R$ such that  $u|_{\p D}\equiv 0$. Let $H^2=W^2_2(D)\cap
H^1$ be the space equipped with the norm of $W_2^2(D)$. The spaces
$H^k$ and $W_2^k(D)$ are called  Sobolev spaces; they are Hilbert
spaces, and $H^k$ is a closed subspace of $W_2^k(D)$, $k=0,1,2$.
\par
 Let $H^{-1}$ be the dual space to $H^{1}$, with the
norm $\| \,\cdot\,\| _{H^{-1}}$ such that if $u \in H^{0}$ then
$\| u\|_{ H^{-k}}$ is the supremum of $(u,v)_{H^0}$ over all $v
\in H^0$ such that $\| v\|_{H^1} \le 1 $. $H^{-k}$ is a Hilbert
space.
\par
We denote by $\ell_k$ and $\oo\ell _{k}$ the Borel measure and the
Lebesgue measure in $\R^k$ respectively, and we denote by ${\cal
B}_{k}$ the $\sigma$-algebra of Borel sets in $\R^k$.  We denote
by $ \oo{{\cal B}}_{k}$ the completion of ${\cal B}_{k}$ with
respect to the measure $\ell_k$, or the $\sigma$-algebra of
Lebesgue sets in $\R^k$.
\par
We denote by $\oo{{\cal P}}$  the completion (with respect to the
measure $\oo\ell_1\times\P$) of the $\s$-algebra of subsets of
$[0,T]\times\O$, generated by functions that are progressively
measurable with respect to $\F_t$.
\par
Let $Q_s\defi D\times [s,T]$. For $k=-1,0,1,2$, we  introduce spaces
 \baaa
 X^{k}(s,T)\defi L^{2}\bigl([ s,T ]\times\Omega,
{\oo{\cal P }},\oo\ell_{1}\times\P;  H^{k}\bigr),\quad Z^k_t
\defi L^2\bigl(\Omega,{\cal F}_t,\P;
H^k\bigr),\quad \CC^{k}(s,T)\defi C\bigl([s,T]; Z^k_T\bigr). \eaaa
The spaces $X^k$ and $Z_t^k$  are Hilbert spaces.
\par
Further, we introduce spaces $$ Y^{k}(s,T)\defi X^{k}(s,T)\!\cap
\CC^{k-1}(s,T), \quad k\ge 0, $$ with the norm $ \| u\| _{Y^k(s,T)}
\defi \| u\| _{{X}^k(s,T)} +\| u\| _{\CC^{k-1}(s,T)}. $
\par
For brevity, we will use the notations
 $X^k\defi X^k(0,T)$, $\CC^k\defi \CC^k(0,T)$,
and  $Y^k\defi Y^k(0,T)$.
\par
In addition, we will be using spaces \baaa &&\Z^k_c\defi
L_2(\O,\F_T,\P;C^k(D)),\quad \X^k_c= L^{2}\bigl([ 0,T ]\times\O,\,
\oo{{\cal P} },\oo\ell_{1}\times\P;\; C^k(\oo D)\bigr),\quad k\ge 0,
\\ &&\W^{k}_p \defi L^{\infty}\bigl([0,T ]\times\O, \overline{{\cal
P}},\oo\ell_{1}\times\P;\, W_p^{k}(D)\bigr), \quad k=0,1,\ldots,
\quad 1\le p\le
 +\infty.
\eaaa  The same notations will be used for the spaces of vector and
matrix functions, meaning that all components belong to the
corresponding spaces.  In particular, $\|\cdot\|_{\W^k_p}$ means the
sum of all this norms for all components.
\par We will write $(u,v)_{H^0}$ for $u\in H^{-1}$
and $v\in H^1$, meaning the obvious extension of the bilinear form
from $u\in H^{0}$ and $v\in H^1$. Similarly, we will write
$(\xi,\eta)_{X^0}$ for $\xi\in X^{-1}$ and $\eta\in X^1$.
\begin{proposition} 
\label{propL} Let $\xi\in X^0$,
 let a sequence  $\{\xi_k\}_{k=1}^{+\infty}\subset
L^{\infty}([0,T]\times\O, \ell_1\times\P;\,C(\oo D))$ be such that
all $\xi_k(\cdot,t,\o)$ are progressively measurable with respect to
$\F_t$, and let $\|\xi-\xi_k\|_{X^0}\to 0$. Let $t\in [0,T]$ and
$j\in\{1,\ldots, N\}$ be given.
 Then the sequence of
integrals $\int_0^t\xi_k(x,s,\o)\,dw_j(s)$ converges in $Z_t^0$ as
$k\to\infty$, and its limit depends on $\xi$, but does not depend
on $\{\xi_k\}$.
\end{proposition}
\par
{\it Proof} follows from completeness of  $X^0$ and from the
equality
\begin{eqnarray*}
\E\int_0^t\|\xi_{k}(\cdot,s,\o)-\xi_m(\cdot,s,\o)\|_{H^0}^2\,ds
=\int_D\,dx\,\E\left(\int_0^t\big(\xi_k(x,s,\o)-
\xi_m(x,s,\o)\big)\,dw_j(s)\right)^2.
\end{eqnarray*}
\begin{definition} 
\rm For $\xi\in X^0$, $t\in [0,T]$, and $j\in\{1,\ldots, N\}$, we
define $\int_0^t\xi(x,s,\o)\,dw_j(s)$ as the limit  in $Z_t^0$ as
$k\to\infty$ of a sequence $\int_0^t\xi_k(x,s,\o)\,dw_j(s)$, where
the sequence $\{\xi_k\}$ is such  as in Proposition \ref{propL}.
\end{definition}
 Sometimes we will omit
$\o$.
\section{Forward and backward SPDEs}
In this section, we collect some known fact for SPDEs.
\subsection{Forward SPDEs}
\label{SecC} Let $s\in [0,T)$, $\varphi\in X^{-1}$, $h_i\in X^0$,
and $\Phi\in Z^0_s$. Consider the boundary value problem \be
\label{4.1} \ba d_tu=\left( \A u+ \varphi\right)dt + \sum_{i=1}^N[
B_iu+h_i]dw_i(t), \quad t\ge s,\\ u|_{t=s}=\Phi,\quad\quad
u(x,t,\o)|_{x\in \p D}=0. \ea
 \ee
 Here
 $u=u(x,t,\o)$,
 $(x,t)\in Q$,   $\o\in\O$, and
  \be\label{A}\A v\defi \sum_{i,j=1}^n
b_{ij}(x,t,\o)\frac{\p^2v}{\p x_i \p x_j}(x) +\sum_{i=1}^n
f_{i}(x,t,\o)\frac{\p v}{\p x_i }(x) +\,\lambda(x,t,\o)v(x), \ee
where $b_{ij}, f_i, x_i$ are the components of $b$,$f$, and $x$.
Further, \be\label{B} B_iv\defi\frac{dv}{dx}\,(x)\,\beta_i(x,t,\o)
+\oo \beta_i(x,t,\o)\,v(x),\quad i=1,\ldots ,N. \ee
\par
We assume that the functions $b(x,t,\o):
\R^n\times[0,T]\times\O\to\R^{n\times n}$, $\b_j(x,t,\o):
\R^n\times[0,T]\times\O\to\R^n$, $\oo\b_i(x,t,\o):$
$\R^n\times[0,T]\times\O\to\R$, $f(x,t,\o):
\R^n\times[0,T]\times\O\to\R^n$, $\lambda(x,t,\o):
\R^n\times[0,T]\times\O\to\R$ and  $\varphi (x,t,\o): \R^n\times
[0,T]\times\O\to\R$ are  progressively measurable for any $x\in
\R^n$ with respect to $\F_t$.\par
 To proceed further, we assume that Conditions
\ref{cond3.1.A}-\ref{condB} remain in force throughout this paper.
 \begin{condition} \label{cond3.1.A} The matrix  $b=b^\top$ is
symmetric,  bounded, and progressively measurable with respect to
$\F_t$ for all $x$, and there exists a constant $\d>0$ such that \be
 \label{Main1} y^\top  b
(x,t,\o)\,y-\frac{1}{2}\sum_{i=1}^N |y^\top\b_i(x,t,\o)|^2 \ge
\d|y|^2 \quad\forall\, y\in \R^n,\ (x,t)\in  D\times [0,T],\
\o\in\O. \ee
\end{condition}
\par
Inequality (\ref{Main1}) is called sometimes a coercivity
condition; it means that equation (\ref{4.1}) is {\it
superparabolic}, in terminology of Rozovskii (1990).
\begin{condition}\label{cond3.1.B} The
functions $b(x,t,\o):\R^n \times \R\times\O\to \R^{n \times n}$,
$f(x,t,\o):\R^n \times \R\times\O\to \R^n$, $\lambda (x,t,\o):\R^n
\times \R\times\O\to \R$,  are bounded and differentiable in $x$,
and \baaa \esssup_{(x,t,\o)\in    Q} \biggl[ \Bigl| \frac{\p b}{\p
x}(x,t,\o)\Bigr| + \Bigl| \frac{\p f}{\p x}(x,t,\o)\Bigr| + \Bigl|
\frac{\p \lambda}{\p x}(x,t,\o)\Bigr| \biggl]< +\infty. \eaaa
\end{condition}
\begin{condition}\label{condB}
The functions  $\b_i(x,t,\o)$ and $\oo\b_i(x,t,\o)$ are bounded
and differentiable in $x$, and  $\esssup_{x,t,\o}|\frac{\p
\b_i}{\p x}(x,t,\o)|<+\infty$, $\esssup_{x,t,\o}|\frac{\p\oo
\b_i}{\p x}(x,t,\o)|<+\infty$, $i=1,\ldots ,N$.
\end{condition}
\par
We introduce the set of parameters $$ \ba {\cal P}_1
\defi \biggl( n,\,\, D,\,\, T,\,\, \ \d,\,\,\,\,
\esssup_{x,t,\o}\Bigl[| b(x,t,\o)|+ |f(x,t,\o)|+ \Bigl|\frac{\p
b}{\p x}(x,t,\o)\Bigr|+
\Bigl|\frac{\p f}{\p x}(x,t,\o)\Bigr|\Bigr],\\
 \esssup_{x,t,\o,i}\Bigl[|
\b_i(x,t,\o)|+ |\oo\b_i(x,t,\o)|+ \Bigl|\frac{\p \b_i}{\p
x}(x,t,\o)\Bigr|+ \Bigl|\frac{\p\oo\b_i}{\p x}(x,t,\o)\Bigr|\Bigr]
\biggr). \ea $$
\subsubsection*{The definition of solution} 
\begin{definition} 
\label{defsolltion} \rm Let $h_i\in X^0$ and $\varphi\in X^{-1}$.
We say that   equations (\ref{4.1}) are satisfied for $u\in Y^1$
if
\begin{eqnarray}
&&u(\cdot,t,\o)-u(\cdot,r,\o)\nonumber
\\  &&\hphantom{xxx}= \int_r^t\big(\A u(\cdot,s,\o)+
\varphi(\cdot,s,\o)\big)\,ds+ \sum_{i=1}^N
\int_r^t[B_iu(\cdot,s,\o)+h_i(\cdot,s,\o)]\,dw_i(s)
\label{intur}
\end{eqnarray}
for all $r,t$ such that $0\le r<t\le T$, and this equality is
satisfied as an equality in $Z_T^{-1}$.
\end{definition}
Note that the condition on $\p D$ is satisfied in the following
sense:   $u(\cdot,t,\o)\in H^1$ for a.e. \ $t,\o$. Further,  the
value of  $u(\cdot,t,\o)$ is continuous in $t$ in $Z_T^0$ and
uniquely defined in $Z_T^0$ given $t$, by the definitions of the
space  $Y^1$. The stochastic integrals with $dw_i$ in (\ref{intur})
are defined as elements of $Z_T^0$. For an arbitrary process $u\in
Y^1$, the integral with $ds$ is defined as an element of $Z_T^{-1}$.
However, $u\in Y^1$ presented in Definition \ref{defsolltion} is
such that this integral is equal to an element of $Z_T^{0}$ in the
sense of equality in $Z_T^{-1}$.
\subsubsection*{Existence and regularity for forward SPDEs}
Typically, existence and uniqueness results at different spaces
for linear PDEs are based on so-called  prior estimates, when a
norm of the solution is estimated via a norm of the free term. For
the second order equations, there are two important estimates
based on $L_2$-norm: so-called "the first energy inequality" or
"the first fundamental inequality", and "the first energy
inequality", or "the second fundamental inequality" (Ladyzhenskaya
(1985)). For instance, consider a boundary value problem for the
heat equation \baaa
&&u_t'=u''_{xx}+\varphi,\quad\varphi=f'_x+g,\nonumber
\\&&u|_{t=0}=0, u|_{\p D}=0,\quad (x,t)\in Q=D\times [0,T], \quad
D\subset \R. \label{simple} \eaaa Then the first fundamental
inequality is the estimate $$\|u\|^2_{L_2(Q)}+
\|u'_x\|^2_{L_2(Q)}\le\const (\|f\|^2_{L_2(Q)}+\|g\|^2_{L_2(Q)}).$$
Respectively, the second fundamental inequality is the estimate $$
\|u'_t\|^2_{L_2(Q)}+\|u\|^2_{L_2(Q)}+\|u'_x\|^2_{L_2(Q)}+\|u''_{xx}\|^2_{L_2(Q)}
\le\const \|\varphi\|^2_{L_2(Q)}.$$ The second fundamental
inequality leads to existence theorem in the class of functions $u$
such that $u''_{xx}\in L_2(Q)$. The first fundamental inequality
allows more general free terms but leads to existence theorem in the
class of functions $u$ with generalized derivatives $u''_{xx} \in
H^{-1}$ only.
\par
An analog of the first and the second fundamental inequality for
the forward SPDEs is given by the following two theorems.
\begin{theorem}
\label{lemma1} [Rozovskii (1990), Ch. 3.4.1] Assume that  Conditions
\ref{cond3.1.A}, \ref{cond3.1.B},  and \ref{condB}, are satisfied.
 Then problem (\ref{4.1}) has an unique solution $u$ in the
class $Y^1(s,T)$ for any $\varphi \in X^{-1}(s,T)$,   $\Phi\in
Z_s^0$, $h_i\in X^0(s,T)$, $i=1,\ldots ,N$, and the following analog
of the first fundamental inequality is satisfied:  \be \label{4.2}
\| u \|_{Y^1(s,T)}\le c \left(\| \varphi \|
_{X^{-1}(s,T)}+\|\Phi\|_{Z^0_s}+ \sum_{i=1}^N\|h_i
\|_{X^0(s,T)}\right), \ee where  $c=c({\cal P}_1)$ is a constant
that depends on ${\cal P}_1 $ only.
\end{theorem}
\begin{theorem}\label{Th3.1.1} [Dokuchaev (2005)]
 Assume that
Conditions \ref{cond3.1.A}, \ref{cond3.1.B},  and \ref{condB},
 are satisfied.
In addition, assume that $\b_i(x,t,\o)=0$ for $x\in \p D$,
$i=1,...,N$.  Then problem (\ref{4.1}) has an unique solution $u\in
{Y}^2$ for any $\varphi \in X^0$, $\Phi \in Z_0^1$,  $h_i\in X^1$,
$i=1,\ldots ,N$,
 and the following analog of the second fundamental
inequality is satisfied: \be \| u \|_{{Y}^2}
 \le   c\left(  \|\varphi \|_{X^0} +
\|\Phi \|_{Z_0^1} +  \sum_{i=1}^N\|h_i \|_{X^1}\right),
\label{3.1.3} \ee where  $c=c({\cal P}_1)$ is a constant that
depends on ${\cal P}_1 $ only.
\end{theorem}
\par
Introduce  operators $L(s,T):X^{-1}(s,T)\to Y^1(s,T)$,
$\M_{i}(s,T):X^{0}(s,T)\to Y^1(s,T)$, and $\L(s,T):Z^0_s\to
Y^1(s,T)$, such that
$$u=L(s,T)\varphi+\L(s,T)\Phi+\sum_{i=1}^N\M_{i}(s,T)h_i,$$ where
$u$ is the solution in $Y^1(s,T)$ of   problem (\ref{4.1}). These
operators are linear and continuous; it follows immediately from
Theorem \ref{lemma1}. We will denote by $L$, $\M_i$, and $\L$, the
operators $L(0,T)$, $\M_i(0,T)$, and $\L(0,T)$, correspondingly.
\subsection{Backward  SPDEs}
\label{secBE} Introduce the operators being formally adjoint to the
operators $\A$ and $B_i$:
\begin{eqnarray*}
\A^*v &=&\sum_{i,j=1}^n\frac{\p ^2 }{\p x_i \p x_j}\,
\biggl(b_{ij}(x,t,\o)\,v(x)\biggr)-\, \sum_{i=1}^n \frac{\p }{\p
x_i }\,\big( f_{i}(x,t,\o)\,v(x)\big) +\lambda(x,u,t,\o)\,v(x),
\\
 B_i^*v&=&-\sum_{i=1}^n \frac{\p }{\p x_i }\,\big(\beta_{i}(x,t,\o)\,v(x))+\oo \beta_i(x,t,\o)\,v(x).
\end{eqnarray*}
Consider the boundary value problem in $Q$ \baa   && d_tp+
\Bigl(\A^*p+\sum_{i=1}^NB_i^*\chi_i+\xi\Bigr)\,dt=
\sum_{i=1}^N\chi_i \,dw_i(t),\nonumber
\\ &&p|_{t=T}=\Psi, \quad p(x,t,\o)\,|_{x\in \p D}=0.
\label{4.8}
\eaa
\subsubsection*{The definition of solution}
\begin{definition} 
\label{defsolltion2} \rm  We say that equation (\ref{4.8}) is
satisfied for $p\in Y^1$, $\xi\in X^{-1}$, $\Psi\in Z_T^0$,
$\chi_i\in X^0$ if
\begin{eqnarray} 
\label{intur1} p(\cdot,t)=\Psi+\int_t^T\Biggl(\A^*p(\cdot,s) +
\sum_{i=1}^NB_i^*\chi_i(\cdot,s)+\xi(\cdot,s)\Biggr)\,ds
-\sum_{i=1}^N \int_t^T\chi_i(\cdot,s)\,dw_i(s)
\end{eqnarray}
for any $t\in[0,T]$. The equality here is assumed to be an
equality in the space $Z_T^{-1}$.
\end{definition}
\subsubsection*{Existence and regularity for backward SPDEs}
For $t\in[0,T]$, define operators $\I_t: C([0,T];Z_T^{k})\to
Z^{k}_t$ such that $\I_tu=u(\cdot,t)$. \par The following theorem
gives an analog of the first fundamental inequality for backward
SPDEs. In addition, this theorem establishes duality between forward
and backward equations.
\begin{theorem}  
\label{Th4.2} [Dokuchaev (1992,2010)] For any $\xi \in X^{-1}$ and
$\Psi\in Z_T^0$,  there exists a pair $(p,\chi)$, such that $p\in
Y^1$, $\chi=(\chi_1,\ldots, \chi_N)$, $\chi_i\in X^0$ and
(\ref{4.8}) is satisfied. This pair is uniquely defined, and the
following analog of the first fundamental inequality is satisfied:
\be\label{1FE} \|p\|_{Y^1}+\sum_{i=1}^N\|\chi_i\|_{X^0}\le
c(\|\xi\|_{X^{-1}}+\|\Psi\|_{Z_T^{0}}), \ee where $ c=c({\cal
P}_1)>0$ is a constant that depends on ${\cal P}_1$ only.
Furthermore, the following duality holds between problems
(\ref{4.8}) and (\ref{4.1}):
$$p={\LU}^*\xi+(\I_T\LU)^*\Psi,\quad
\chi_i=\M_i^*\xi+(\I_T\M_i)^*\Psi,\quad
p(\cdot,0)=\L^*\xi+(\I_T\L)^*\Psi,$$ where ${\LU}^*: X^{-1}\to X^1$,
${\M}_i^*: X^{0}\to X^0$, $(\I_T\LU)^*: Z_0^{0}\to X^1$,
$(\I_T\M_i)^*: Z_0^{0}\to X^0$, and $(\I_T\L)^*: Z_T^{0}\to Z_0^0$,
 are the operators that are adjoint  to the operators ${\LU}:
X^{-1}\to X^1$, $\M_i : X^{0}\to X^1$, $\I_T\M_i: X^{-1}\to
Z_T^0$,  $\I_T\M_i: X^0\to Z_T^{0}$, and $\,\I_T\L: Z_0^0\to
Z_T^{0}$, respectively.
\end{theorem}
\par
We will need an analog of the  second fundamental inequality as
well. \par Starting from now, we assume that the following addition
conditions are satisfied.
\begin{condition}\label{cond5}
There exist functions $\w f(x,t,\o): \R^n\times\R_+\times\O\to
\R^n$, $\w \lambda (x,t,\o): \R^n\times\R_+\times\O\to \R$, and $\w
\b_i(x,t,\o): \R^n\times\R_+\times\O\to \R$, such that \baaa
 \esssup_{x,t,\o}\Bigl(|\w f(x,t,\o)|+|\w \lambda (x,t,\o)|+|\w \b_i(x,t,\o)|\Bigr)<+\infty,
 \eaaa
and
 \baaa
\A^*p &=& \sum_{i,j=1}^n b_{ij}(x,t,\o)\,\frac{\p ^2 p}{\p x_i \p
x_j}\,(x) + \sum_{i=1}^n \hat f_{i}(x,t,\o)\,\frac{\p p}{\p x_i
}\,(x)- \w \lambda(x,t,\o)\,p(x),
\nonumber\\
B_i^*p &=&\frac{d p}{d x}(x)\,\b_i(x,t,\o)+\w\beta_i(x,t,\o)\,p(x).
\label{B*}\eaaa
\end{condition}
Clearly, this condition is satisfied if  the function
$b(x,t,\o):\R^n \times \R\times\O\to \R^{n \times n}$ is twice
differentiable in $x$, and \baaa \esssup_{\o}\sup_{(x,t)\in Q}
\Bigl| \frac{\p^2 b}{\p x_k\p x_m}(x,t,\o)\Bigr|< +\infty. \eaaa
 \par
For an integer $M>0$, let $\Theta_b(M)$ denotes  the class of all
matrix functions  $b$  such that all conditions imposed in Section
\ref{SecC} are satisfied, and there exists a set $\{T_i\}_{i=0}^M$
such that $0=T_0<T_1<\cdots <T_M=T$ and that the function
$b(x,t,\o)=b(x,\o)$  does not depend on $t$ for $t\in[T_i,T_{i+1})$.
(it follows from the assumptions that $b(x,t,\cdot)$ is
$\F_{T_i}$-measurable for all $x\in D$, $t\in [T_i,T_{i+1})$).
\par Let $\Theta_b\defi
\cup_{M>0} \Theta_b(M)$.
\par
Let $\oo\Theta_b$ denotes the class of function $b$ from such that
all conditions imposed in Section \ref{SecC} are satisfied,  and
there exists  and a sequence
$\{b^{(i)}\}_{i=1}^{+\infty}\subset\Theta_b$ such that
$\|b-b^{(i)}\|_{\W_{\infty}^1}\to 0$ as  $i\to +\infty$. (Remind
that the assumptions on $b$ are such that $b\in\W_{\infty}^1$).
\par
\begin{condition} \label{condA}  The
matrix $b$ belongs to $\oo\Theta_b$,
 and  there exists a constant $\d_1>0$ such that
\baa
 \label{Main1'} \sum_{i=1}^Ny_i^\top  b
(x,t,\o)\,y_i-\frac{1}{2}\left(\sum_{i=1}^Ny_i^\top\b_i(x,t,\o)\right)^2
\ge \d_1\sum_{i=1}^N|y_i|^2\nonumber\\ \quad \forall\,
\{y_i\}_{i=1}^N\subset \R^n,\ (x,t)\in D\times [0,T],\ \o\in\O. \eaa
\end{condition}
\begin{remark} If Condition \ref{condA} holds, then
Condition \ref{cond3.1.A} holds. If $n=1$ and Condition
\ref{cond3.1.A} holds, then the estimate in Condition \ref{condA}
also holds. If $n>1$, then it can happen that Condition
\ref{cond3.1.A} holds, but the estimate in Condition \ref{condA}
does not hold. For instance, assume that $n=2$, $N=2$, $\b_1\equiv
(1,0)^\top$, $\b_2=(0,1)^\top$, $b\equiv
\frac{1}{2}(\b_1\b_1^\top+\b_2\b_2^\top)+0.01I_2=0.51I_2$, where
$I_2$ is the unit matrix in $\R^{2\times 2}$. Obviously, Condition
\ref{cond3.1.A} holds and $b\in\oo\Theta_b$,. On the other hand,
Condition \ref{condA} does not hold for this $b$;  to see this, it
suffices to take $y_1=\b_1$ and $y_2=\b_2$.
\end{remark}
\begin{remark}
 Condition \ref{condA} is satisfied for matrices $b\in\oo\Theta_b$
 if either $n=1$ or
 there exists $N_0\in\{1,...,N\}$ such that
 $\b_i\equiv 0$ for $i>N_0$, and  there exists a constant $\d_2>0$ such that
\be
 \label{Main1''} y^\top  b
(x,t,\o)\,y-\frac{N_0}{2}|y^\top\b_i(x,t,\o)|^2 \ge \d_2|y|^2
\quad\forall\, y\in \R^n,\ (x,t)\in  D\times [0,T],\ \o\in\O,\
i=1,...,N_0. \ee In particular, it is satisfied if Condition
\ref{cond3.1.A} holds and $N_0=1$.
\end{remark}
To proceed further, we assume that Conditions \ref{cond5}-
\ref{condA} remain in force starting from here and up to the end of
this paper, as well as
the previously formulated conditions.
\par
Let ${\cal P}\defi ({\cal P}_1,\d_1)$.
\par
We will be using the following analog of the second fundamental
inequality for backward SPDEs.
\begin{theorem} 
\label{Th2FE}  [Dokuchaev (2006)] For any $\xi \in X^{0}$ and
$\Psi\in Z_T^1$,  there exists a pair $(p,\chi)$, such that $p\in
Y^2$, $\chi=(\chi_1,\ldots, \chi_N)$, $\chi_i\in X^1$ and
(\ref{4.8}) is satisfied. This pair is uniquely defined, and
$$p={\LU}^*\xi+(\I_T\LU)^*\Psi,\quad
\chi_i=\M_i^*\xi+(\I_T\M_i)^*\Psi.$$ The operators $L^*: X^0\to
Y^2$,  $(\I_TL)^*:Z_T^1\to Y^2$, and $\M_i^*:X^0\to X^1$,
$(\I_T\M_i)^*:Z_T^1\to X^1$, are continuous. More precisely, the
following analog of  the second fundamental inequality holds:
\be\label{2FE} \|p\|_{Y^2}+\sum_{i=1}^N\|\chi_i\|_{X^1}\le
c(\|\xi\|_{X^{0}}+\|\Psi\|_{Z_T^{1}}), \ee where $ c>0$ is a
constant that depends only on ${\cal P}$.
\end{theorem}
\subsubsection*{Semi-group property  for
backward equations} It is known that the dynamic of forward
parabolic Ito equation has semi-group property (or causality
property): if $u=L\varphi+\L_0\Phi$, where $\varphi\in X^{-1}$,
$\Phi\in Z_0^0$, then \be\label{casF}
u|_{t\in[\t,s]}=(L\varphi+\L_0\Phi)|_{t\in[\t,s]}=L(\t,s)\varphi+\L_{\t}(\t,s)u(\cdot,\t).\ee
We will need  a similar property for the backward equations.
\begin{theorem}\label{lemmaMark} (Semi-group property for backward equations)   [Dokuchaev (2010)].
Let $0\le \t< s<T$, and let $p=L^*\xi$, $\chi_i=\M_i\xi$ where
$\xi\in X^{-1}$ and $\Psi\in Z_T^0$. Then \baa
&&p|_{t\in[\t,s]}=L(\t,s)^*\xi|_{t\in[\t,s]}+(\I_sL(\t,s))^*p(\cdot,s),\label{p-cas}
\\
&&p(\cdot,\t)=(\I_s\L_{\t}(\t,s))^*p(\cdot,s)+\L_{\t}(\t,s)^*\xi,\label{p-cas0}
\\
&&\chi_i|_{t\in[\t,s]}=\M_i(\t,s)^*\xi|_{t\in[\t,s]}+(\I_s\M_i(\t,s))^*p(\cdot,s),
\quad k=1,...,N. \label{chi-cas} \eaa
\end{theorem}
\subsubsection{Some additional regularity}
Theorem \ref{Th2FE} requires that $\Psi\in Z_2^1$. We will need a
modification of this theorem that allows $\Psi\in Z_2^0$.
\begin{theorem}
\label{Th3FE} Let the assumptions of Theorem \ref{Th2FE} be
satisfied. Let $\xi \in X^{0}$ and $\Psi\in Z_T^0$. Let
$$p={\LU}^*\xi+(\I_T\LU)^*\Psi,\quad
\chi_i=\M_i^*\xi+(\I_T\M_i)^*\Psi.$$ Let $\e\in(0,T)$ be given. Then
\be\label{3FE}
\|p\|_{Y^2(0,T-\e)}+\sum_{i=1}^N\|\chi_i\|_{X^1(0,T-\e)}\le
\frac{c}{\sqrt{\e}}(\|\xi\|_{X^{0}}+\|\Psi\|_{Z_T^{0}}), \ee where
$c= c({\cal P})>0$ is a constant that depends only on and ${\cal
P}$.
\end{theorem}
\par {\it Proof.}
By Theorem \ref{Th4.2} and Theorem \ref{lemmaMark}, it follows that
\be\label{4FE}
\|p\|_{Y^2(0,T-\e)}+\sum_{i=1}^N\|\chi_i\|_{X^1(0,T-\e)}\le
c_1(\|\xi\|_{X^{0}(0,T-\e)}+\|p(\cdot,T-\e)\|_{Z_T^{1}}), \ee where
$ c_1=c_1({\cal P})>0$ is a constant that depends only on ${\cal
P}$. (Note that the same constant $c$ can be used for all $\e$,
since Theorem \ref{4.2} holds for $T$ replaced by $T-\e$ with any
$\e\in[0,T)$). In addition, it follows from Theorem \ref{Th2FE} that
$$ \inf_{s\in [T-\e,T]}\|p(\cdot,s)\|^2_{Z_T^{1}}\le
\frac{1}{\e}\int_{T-\e}^T\|p(\cdot,t)\|^2_{Z^1_T}dt\le
\frac{c_2}{\e}(\|\xi\|_{X^{0}}^2+\|\Psi\|_{Z_T^{0}}^2), $$ where $
c_2=c_2({\cal P})>0$ is a constant that depends only on ${\cal P}$.
This completes the proof. $\Box$
\section{The main result: the representation theorem}
Let functions ${\ww\b_i: Q\times \O \to \R^n}$, $i=1,\ldots, M$, be
such that $$ 2b(x,t,\o)=\sum_{i=1}^N\b_i(x,t,\o)\,\b_i(x,t,\o)^\top
+\sum_{j=1}^M\,\ww\b_j(x,t,\o)\,\ww\b_j(x,t,\o)^\top, $$ and $\ww
\b_i$ has the similar properties as $\b_i$. (Note that, by Condition
\ref{cond3.1.A}, $2b>\sum_{i=1}^N\b_i\b_i^\top$).
\par
 Let
$\ww w(t)=(\ww w_1(t),\ldots, \ww w_M(t))$ be a new Wiener process
independent on $w(t)$.
\par
Let  $(x,s)\in \oo D\in[0,T]$  be given. Consider the following Ito
equation
\begin{eqnarray}
\label{yxs} &&dy(t) =\ww
f(y(t),t)\,dt+\sum_{i=1}^N\b_i(y(t),t)\,dw_i(t) +\sum_{j=1}^M\ww
\b_j(y(t),t)\,d \ww w_j(t),
\nonumber\\ [-6pt] &&y(s)=x,
\end{eqnarray}
where  $\ww f\defi\hat f-\sum_{i=1}^N\w \b_i\b_i$.
\par
Let  $y(t)=y^{x,s}(t)$ be the solution of (\ref{yxs}). \par Set
$\tau^{x,s}\defi\min\,\{t\le T: y^{x,s}(t)\notin D\}$. For $t\ge s$,
set
\begin{eqnarray*}
\g^{x,s}(t) \defi\exp\biggl[-\int_s^t \w \lambda(y^{x,s}(t),t)\,dt
+\sum_{i=1}^N\int_s^t\w \b_i(y^{x,s}(s),s)\,dw_i(s)
-\sum_{i=1}^N\frac{1}{2}\,\int_s^t\w\b_i(y^{x,s}(s),s)^2\,ds\biggr].
\end{eqnarray*}
\begin{theorem}\label{Threpr}
Let $b\in \X_c^3$, $\w f\in\X_c^2$, $\w\lambda\in\X^1_c$, $
\b_i\in\X_c^3$ and $ \w \b_i\in\X_c^2$.
 Let $(p,\chi_1,...,\chi_N)$ be the solution of (\ref{4.8}), where functions $\xi:Q\times \O\to\R$
and $\Psi:D\times \O$ are such that $\xi$ is $({\cal B}_{n+1}\otimes
\F,{\cal B}_1)$-measurable, $\Psi$ is $({\cal B}_{n}\otimes \F,{\cal
B}_1)$-measurable, $\xi\in X^0$ and $\Psi\in Z_T^0$.
 Then for any
$s\in[0,T)$,  \baa\label{repr}
p(x,s,\o)=\E\Big\{\g^{x,s}(T)\Psi(y^{x,s}(T))\Ind_{\{T\ge
\tau^{x,s}\}}\,|\,\F_s\Big\}+\E\biggl\{\int_s^{\tau^{x,s}}
\g^{x,s}(t)\,\xi(y^{x,s}(t),t,\o)\,dt\,\bigr|\,\F_s\biggr\}\\
\hbox{for a.e.}\ x,\o. \nonumber\eaa
\end{theorem}
\par
Remind that the solution   $(p,\chi_1,...,\chi_N)$  of (\ref{4.8})
can be represented as  \be
 p= L^*\xi+(\I_T L)^*\Psi,\qquad \chi_i=\M_i^*\xi+(\I_T \M_i)^*\Psi,\quad i=1,...,N.
 \label{pintheorem}\ee
\section{Proof of  Theorem \ref{Threpr}}
 Let us proof first the following lemma.
\begin{lemma}\label{lemmatemp}  Theorem \ref{Threpr} holds
 even without  Condition \ref{cond3.1.A} for the case when $\xi\in \X^0_c$,
$\Psi\in\Z^0_c\cap Z^1_0$, $p\in\X^2_c$,
 $p(\cdot,T)\in Z_c^0$, $\chi_i\in\X_c^1$,
 where $(p,\chi_1,....,\chi_N)$ is the solution of (4.8).
\end{lemma}
\par
{\it Proof of Lemma \ref{lemmatemp}}.
 Let $(x,s)$ be given, and let
$y(t)=y^{x,s}(t)$ and  $\g(t)=\g^{x,s}(t)$. We have that $$ d_t
p={\cal J}(p)\,dt + \sum_{i=1}^N\chi_i\,dw_i(t), $$ where
$$ {\cal J}(p)\defi -\A^*p-\sum_{i=1}^NB_i^*\chi_i-\xi. $$ Let
$\psi(t)\defi p(y(t,\o),t,\o)$.
\par
By the Ito-Ventssel formula (see, e.g., Rozovskii (1990), Chapter 1
), \baaa d\psi(t) = h(y(t),t)dt+ \,\sum_{i=1}^N\chi_i
(y(t),t)\,dw_i(t)+\,\sum_{i=1}^N \bigg(\frac{\p p}{\p
x}\b_i\bigg)\,(y(t),t)\,d w_i(t)\\  + \sum_{i=1}^M\bigg(\frac{\p
p}{\p x}\ww\b_i\bigg)\,(y(t),t)\,d \ww w_i(t),\eaaa where \baaa
h=h(y(t),t) ={\cal J}(p)+\A^*p+\w\lambda p -\frac{\p p}{\p x}\,
\sum_{i=1}^N\w\b_i\b_i + \sum_{i=1}^N \frac{\p \chi_i}{\p x}\b_i.
\eaaa By (\ref{B*}), it can be rewritten as
 \baaa h=-\xi +\w\lambda p -\frac{\p p}{\p x}\,\sum_{i=1}^N\w\b_i\b_i
-\sum_{i=1}^N\w \b_i\chi_i. \eaaa Let $\w\psi(t)\defi
\psi(t)\,\g(t)$, $t\ge s$. We have that \baaa
d\g(t)=\g(t)\Bigl(-\w\lambda
dt+\sum_{i=1}^N\w\b_i(t)\,dw_i(t)\Bigr).\eaaa
 Using Ito formula, we derive
that
$$ d\w\psi(t)=-\g(t)\,\xi(y(t),t,\o) +\sum_{i=1}^N\mu_i(t)\,dw_i(t)
+\sum_{i=1}^M\ww\mu_i(t)\,d\ww w_i(t), $$ where  $\mu_i(\cdot)$ and
$\ww \mu_i(\cdot)$ are some $L_2$--integrable  processes such that
$\mu_i(t)$ and $\ww \mu_i(t)$ are independent from $w_j(r)-w_j(t)$
and  $\ww w_k(r)-\ww w_k(t)$ for all $r>t$, $j,k$. It follows that
\begin{eqnarray*}
\E\Big\{\g(T)\Psi(y(T))\Ind_{\{T\le
\tau^{x,s}\}}\,|\,\F_s\Big\}-p(x,s,\o)
 &=&\E\Big\{\left(p(y(\tau^{x,s}),\tau^{x,s},\o)-p(x,s,\o)\right)
\,|\,\F_s\Big\}
\\
&=& -\E\biggl\{\int_s^{\tau^{x,s}}
\g(t)\,\xi(y(t),t,\o)\,dt\,\bigr|\,\F_s\biggr\}.
\end{eqnarray*}
Then (\ref{repr}) follows. This completes the proof of Lemma
\ref{lemmatemp}. $\Box$
\par
 \vspace{3mm}
Let us continue the proof of Theorem \ref{Threpr}, and  let us
assume first  that the functions $\xi$ and $\Psi$ are bounded. In
addition, we assume for the case when $D=\R^n$  that there exists a
bounded domain $\w D\subset \R^n$ such that $\xi(x,t,\o)=0$ and
$\Psi(x,\o)=0$ for all $x\notin \w D$ for all $t,\o$.
\par
For functions $h\in X^0$, we introduce some transforms  $h_m$,
$m=1,2,..$.
\begin{itemize}
\item[(a)] Let $D\neq \R^n$. In this case, we introduce an orthonormal basis
$\{v_k\}_{k=1}^{\infty}$ in $L_2(D)$
 consisting of the eigenfunctions  for the eigenvalue  problem
 \baa
\Delta v-v=-\lambda v,\quad   v|_{\p D}=0.
 \label{eigen}\eaa
 Here $\Delta$ is the Laplacian. It is known that
 $v\in C^2(\oo D)\cap H^2$ (see, e.g., Theorem III.3.2 from Ladyzhenskaya and Ural'tseva (1968)).
For a function $h\in X^0$, we denote by $h_m$ the function $h_m\in
X^0$ such that $h_m(\cdot,t,\o)$ is the projection of
$h(\cdot,t,\o)$ on the subspace of $L_2(D)$ generated as the span of
the functions $\{v_k\}_{k=1}^m$.
\item[(b)] Let $D=\R^n$. In this case,  for a function $h\in X^0$, we
denote by $h_m$ the function
$(h)_m(y,t,\o)\defi\int_{\R^n}h(x,t,\o)J^{(m)}(y-x)dx$ the
corresponding Sobolev transform. Here $J(x):\R^n\to\R$ is the
Sobolev kernel: $J(x)=0$ if $|x|\ge 1$, and
$J(x)=\exp\{-|x|/(1-|x|)\}$ if $|x|<1$, and $J^{(m)}(x)=\kappa_n
m^{n}J(mx)$, where $\kappa_n>0$ is such that
$\int_{\R^n}J^{(m)}(x)dx=1$.
\end{itemize}
\par
The transforms $h_m$ has the following properties: \baa && h_m\in
\X_c^2 \quad\forall h\in X,\nonumber\\&&
(h_m,g)_{X^0}=(h,g_m)_{X^0},\quad\forall h,g\in X^0,\nonumber\\ &&
\|h_m\|_{X^1}\le c \|h\|_{X^1}\quad \forall h\in
X^1,\label{h1h1}\eaa for a constant $c>0$ that does not depend on
$h$. The first two properties are obvious.  For the case when
$D=\R^n$, the last property  follows from the known properties of
the Sobolev transform. It suffices to prove the last property for
the case when $D\neq \R^n$. Let $D\neq \R^n$. For any $V\in H^0$, we
have that $V=\sum_{k=1}^{\infty}c_kv_k$, where $c_k=(V,v_k)_{H^0}$,
meaning the convergence of the series in $H^0$. Hence
\baaa\|V_m\|_{H^1}^2=(V_m,V_m-\Delta
V_m)_{H^0}=\Bigl(\sum_{k=1}^{m}c_kv_k,
\sum_{k=1}^{m}c_kv_k-\Delta\sum_{k=1}^{m}c_kv_k\Bigr)_{H^0}\\=
\Bigl(\sum_{k=1}^{m}c_kv_k,
\sum_{k=1}^{m}c_kv_k+\sum_{k=1}^{m}\lambda_kc_kv_k\Bigr)_{H^0}=
\sum_{k=1}^{m}|c_k|^2(1+\lambda_k)\le \|V\|_{H^1}^2.
 \eaaa
 Here $\lambda_k$ are the eigenvalues of problem
 (\ref{eigen}) that correspond to the eigenfunctions $v_k$.
 It follows that (\ref{h1h1}) holds for $D\neq \R^n$.  Therefore, (\ref{h1h1}) holds.
\par
Let $(p,\chi_1,...,\chi_N)\in Y^1\times (X^0)^N$ be such that
 $p=L^*\xi+(\I_TL_i)^*\Psi$ and $\chi_i=\M_i^*\xi+(\I_T\M_i)^*\Psi$.
 By Theorem \ref{Th3FE}, it follows that
 \baa
 (p,\chi_1,...,\chi_N)|_{t\in[0,T-\e)}\in Y^2(0,T-\e)\times (X^1(T-\e))^N\quad \forall\e>0.
 \label{X2e}
 \eaa In particular,
 it follows
that  $\frac {\p^k p}{\p x_k^2}(\cdot,t)_{t\in[0,T-\e)}$, $k=0,1,2$,
and $\frac {\p \chi_i}{\p x_i}_{t\in[0,T-\e)}$ belong to
$X^{0}(0,T-\e)$.
\par
 We have that \baaa &&d_t p_m+[(\A^*
p)_m+\xi_m+\sum_{i=1}^N(B^*_{i}\chi_i)_m]dt=\sum_{i=1}^N\chi_{im}dw_i(t),\\
&&p_m(x,T,\o)=\Psi_m(x,\o),\quad p_m|_{x\in D}=0. \eaaa It can be
rewritten as
 \baaa &&d_t
p_m+[\A^*
p_m+\w\xi^{(m)}+\sum_{i=1}^NB^*_{i}\chi_{im}]dt=\sum_{i=1}^N\chi_{im}dw_i(t),\\
&&p_m(x,T,\o)=\Psi_m(x,\o),\quad p_m|_{x\in \p D}=0. \eaaa Here $$
\w\xi^{(m)}\defi \xi_m+\eta^{(m)},\quad \eta^{(m)}\defi (\A^*
p)_m-\A^* p_m+ \sum_{i=1}^N(B^*_{i} \chi_i)_m- \sum_{i=1}^NB^*_{i}
\chi_{im}.$$
\par
 Let
us show that \be \Psi_m\to\Psi \quad\hbox{in}\ \
Z_T^0\quad\hbox{as}\quad m\to+\infty.\label{xipsi}\ee \par Clearly,
$\Psi_m(\cdot,\o)\to \Psi(\cdot,\o)$ in $L_2(D)$ a.s. In addition,
we have that
$\|\Psi_{m}(\cdot,\o)\|_{L_2(D)}\le\|\Psi(\cdot,\o)\|_{L_2(D)}$.
Hence $\|\Psi_{m}(,\o)-\Psi(\cdot,\o)\|_{L_2(D)}\le
2\|\Psi(,\o)\|_{L_2(D)}$. By the Lebesgue's Dominated Convergence
Theorem, it follows that (\ref{xipsi}) holds. Similarly, we obtain
that \be \xi_m\to\xi \quad\hbox{in}\ \ X^0\quad\hbox{as}\quad
m\to+\infty.\label{xixi}\ee \par  Again, we have
$\xi_m(\cdot,t,\o)\to \xi(\cdot,t,\o)$ in $L_2(D)$ for a.e.
$(t,\o)$. In addition, we have that
$\|\xi_{m}(\cdot,t,\o)\|_{L_2(D)}\le\|\xi(\cdot,t,\o)\|_{L_2(D)}$
and $\|\xi_{m}(,t,\o)-\xi(\cdot,t,\o)\|_{L_2(D)}\le
2\|\xi(,t,\o)\|_{L_2(D)}$. By the Lebesgue's Dominated Convergence
Theorem again, it follows that (\ref{xixi}) holds.
\par
Let us show that \be \w\xi^{(m)}\defi
\xi_m+\eta^{(m)}\to\xi\quad\hbox{weakly in}\ \ X^{-1}\quad
\hbox{as}\ m\to +\infty. \label{wxixiW} \ee By (\ref{xixi}), it
suffices to show that \baa \eta^{(m)}\to 0\quad\hbox{weakly in}\ \
X^{-1}\quad \hbox{as}\ m\to 0. \label{etaetaW}\eaa
\par First, let us show that there exists a constant $c>0$ such that
\be \|\eta^{(m)}\|_{X^{-1}}\le c\quad  \forall m>0. \label{ogr}\ee
 \par
 By Theorem
\ref{Th4.2}, it follows that $\|p\|_{X^1}\le \const$. Hence
$\|p_m\|_{X^{1}}\le \const$.  Hence \baa\|\A^*p_m\|_{X^{-1}}\le
\const.\label{ineq1}\eaa
\par
Further, let $B(X)$ denote the unit ball in a linear normed space
$X$, i.e., $B(X)\defi\{x\in X:\ \|x\|_{X}\le 1\}$. We have that \baa
\|(\A^*p)_m\|_{X^{-1}}=\sup_{y\in
B(X^1)}(y,(\A^*p)_m)_{X^0}=\sup_{y\in B(X^1)}(y_m,\A^*p)_{X^0}\le
\sup_{y\in B(X^1)}\|\A y_m\|_{X^{-1}}\|p\|_{X^1}\nonumber\\ \le c_1
\sup_{y\in B(X^1)}\| y_m\|_{X^{1}}\|p\|_{X^1} \le c_2 \sup_{y\in
B(X^1)}\| y\|_{X^{1}}\|p\|_{X^1}\le c_3.\hphantom{xxxx}
 \label{ineq2}\eaa
Here $c_k$, $k=1,2,3$, are some constant that are independent from
$m$.
\par
Similarly, we have that, by Theorem \ref{Th4.2},
$\|\chi_i\|_{X^0}\le \const$. Hence $\|\chi_{im}\|_{X^{0}}\le
\const$. Hence \baa\|B_i^*\chi_{im}\|_{X^{-1}}\le
\const.\label{ineq3}\eaa
\par
Further, we have that \baa \|(B_i^*\chi_i)_m\|_{X^{-1}}=\sup_{y\in
B(X^1)}(y,(B_i^*\chi_i)_m)_{X^0}=\sup_{y\in
B(X^1)}(y_m,B_i^*\chi_i)_{X^0}\le \sup_{y\in B(X^1)}\|B_i
y_m\|_{X^{0}}\|\chi_i\|_{X^0}\nonumber\\ \le c_1 \sup_{y\in
B(X^1)}\| y_m\|_{X^{1}}\|\chi_i\|_{X^0} \le c_2 \sup_{y\in B(X^1)}\|
y\|_{X^{1}}\|\chi_i\|_{X^0}\le c_3.\hphantom{xxxx}
 \label{ineq4}\eaa
Here $c_k$, $k=1,2,3$, are some constant that are independent from
$m$. Combining (\ref{ineq1})-(\ref{ineq4}), we obtain (\ref{ogr}).
\par
Let $q=q(x,t,\o)$ denote any one of the functions $p$, $\chi_i$, $\p
p/\p x_k$, $\p^2p/\p x_k\p x_m$, $\p\chi_i/\p x_k$, $k,m=1,...,n$,
$i=1,...,N$, $t<T$. Let $\a$ denote the coefficient such that $\a q$
is presented in the expressions $\A^*p$ or $B_i^*\chi_i$.
 \par
For $\t\in[0,T)$, let $X^1(\t)\defi\{h\in X^1: h(\cdot,t)\equiv
0,\quad t\in [\t,T]\}$.
\par
Let $\t\in[0,T]$ and let $h\in X^1(\t)$. It can be shown similarly
to (\ref{xixi}) that  \baaa \a h_m-(\a h)_m\to 0\quad\hbox{in}\quad
X^0\quad\hbox{as} \quad m\to +\infty.\eaaa It follows that \baaa
  &&((\a q)_m-\a
q_m,h)_{X^0}= ((\a q)_m-\a q_m,h)_{X^0(0,\t)}=(q,\a h_m-(\a
h)_m)_{X^0(0,\t)}\to 0\quad\hbox{as} \ m\to \infty.\eaaa We have
that $\eta^{(m)}$ is a sum of different terms expressed as $(\a
q)_m-\a q_m$. Hence $$
  (\eta^{(m)},h)_{X^0(0,\t)}=(\eta^{(m)},h)_{X^0}\to 0\quad\hbox{as} \ m\to +\infty\qquad \forall h\in X^1(\t).
  $$
  Clearly, the set $\cup_{\t\in [0,T)}X^1(\t)$ is dense in
$X^1$.
 By (\ref{ogr}), it follows that
(\ref{etaetaW}) holds. This
  completes the proof of (\ref{wxixiW}).
\par
Let $s\in[0,T)$ be given.
\par
 By (\ref{xixi}), (\ref{wxixiW}),  and
Theorem \ref{Th4.2},  it follows that \baa \ww p_m(\cdot,s)\to
p(\cdot,s)\quad \hbox{weakly in}\ \ Z_T^0\quad\hbox{as}\ \ m\to
0.\label{ppweak}\eaa
 By Mazur's Theorem  (Theorem 5.1.2 from
Yosida (1995)), there exists a sequence of integer numbers $k=k_i\to
+\infty$ such that there exists sets of real numbers
$\{a_{mk}\}_{m=1}^k\subset[0,1]$ such that $\sum_{m=1}^ka_{mk}=1$
and that \baa &&\ww\xi^{(k)}\defi \sum_{m=1}^ka_{mk}\w\xi^{(m)}\to
\xi\quad \hbox{in}\quad X^{-1}\quad \hbox{as}\quad k=k_i\to
+\infty,\nonumber\\&&\ww\Psi^{(k)}\defi \sum_{m=1}^ka_{mk}\Psi_m\to
\Psi\quad \hbox{in}\quad Z^{0}_T\quad \hbox{as}\quad k=k_i\to
+\infty, \nonumber\\&&\ww p^{(k)}(\cdot,s)\defi
\sum_{m=1}^ka_{mk}p_m(\cdot,s)\to p(\cdot,s)\quad \hbox{in}\quad
Z^{0}_T\quad \hbox{as}\quad k=k_i\to +\infty. \label{Mazur} \eaa
Here $\ww p^{(k)}\defi \sum_{m=1}^ka_{mk}p_m$.
\par
By Lemma \ref{lemmatemp}, it follows that, for all $s$ and for a.e.
$x,\o$, \baaa p_m(x,s,\o)=
\E\Big\{\g^{x,s}(T)\Psi_m(y^{x,s}(T))\Ind_{\{T\ge
\tau^{x,s}\}}\,|\,\F_s\Big\}+\E\Bigl\{\int_s^{\tau^{x,s}}
\g^{x,s}(t)\,\w\xi^{(m)}(y^{x,s}(t),t,\o)\,dt\,\Bigr|\,\F_s\Bigr\},
\eaaa and \baaa \ww p^{(k)}(x,s,\o)=
\E\Big\{\g^{x,s}(T)\Psi^{(k)}(y^{x,s}(T))\Ind_{\{T\ge
\tau^{x,s}\}}\,|\,\F_s\Big\}+\E\Bigl\{\int_s^{\tau^{x,s}}
\g^{x,s}(t)\,\ww\xi^{(k)}(y^{x,s}(t),t,\o)\,dt\,\Bigr|\,\F_s\Bigr\},
\eaaa
 By
the assumptions about the boundedness and  the type of
measurability of the functions $\xi:Q\times\O\to\R$ and
$\Psi:D\times\O\to\R$, it follows that $\g^{x,s}(T)\Psi(y(T))
\Ind_{\{T\le \tau^{x,s}\}}$ and $\int_s^{\tau^{x,s}}
\g^{x,s}(t)\,\xi(y(t),t,\o)\,dt$ are bounded random variables. Let
\be \ww p(x,s,\o)\defi
\E\Big\{\g^{x,s}(T)\Psi(y^{x,s}(T))\Ind_{\{T\ge
\tau^{x,s}\}}\,|\,\F_s\Big\}+\E\Bigl\{\int_s^{\tau^{x,s}}
\g^{x,s}(t)\,\xi(y^{x,s}(t),t,\o)\,dt\,\Bigr|\,\F_s\Bigr\}. \ee
Clearly, $\ww p(\cdot,s)\in Z^0_T$.
\par
Let us show that, for a given $s$, \baa \ww p^{(k)}(\cdot,s)\to\ww
p(\cdot,s)\quad \hbox{weakly in}\ \ Z_T^0\quad\hbox{as}\ \ m\to
\infty\label{wpp}.\eaa By (\ref{Mazur}),  property (\ref{wpp})
implies that $p=\w p$. Therefore, if we prove (\ref{wpp}) then
Theorem \ref{Threpr} will be proved for the case when the functions
$\Psi$ and $\xi$ are bounded and finitely (in $x$) supported.
\par
Let us prove (\ref{wpp}).\par
 Without a loss
of generality, we assume that $\Psi(x,\o)=0$,  $\Psi_m(x,\o)=0$,
$\xi(x,t,\o)=0$, $\w\xi^{(m)}(x,t,\o)=0$ for all $x\notin\oo D$. It
follows that  $\Psi^{(k)}(x,\o)=0$ and  $\ww\xi^{(k)}(x,t,\o)=0$ for
all $x\notin\oo D$.
 Let $\rho\in Z_s^0$. We have that \baaa
|(\ww p^{(k)}(\cdot,s)-\ww p(\cdot,s),\rho)_{Z_T^0}|\le
\E\int_{D}\rho(x)\E\Big\{\g^{x,s}(T)|\Psi^{(k)}(y^{x,s}(T))-\Psi(y^{x,s}(T))|\Ind_{\{T\ge
\tau^{x,s}\}}\,|\,\F_s\Big\}dx\\
+\E\int_{D}\rho(x)\E\Bigl\{\int_s^{\tau^{x,s}}
\g^{x,s}(t)\,|\ww\xi^{(k)}(y^{x,s}(t),t,\o)-
\xi(y^{x,s}(t),t,\o)|\,dt\,\Bigr|\,\F_s\Bigr\}dx
\\
\le \E\int_{D}\rho(x)
\E\Big\{\g^{x,s}(T)|\Psi^{(k)}(y^{x,s}(T))-\Psi(y^{x,s}(T))|\,|\,\F_s\Big\}dx
\\+\E\int_{D}\rho(x)\E\Bigl\{\int_s^T
\g(t)\,|\ww\xi^{(k)}(y^{x,s}(t),t,\o)-\xi(y^{x,s}(t),t,\o)|\,dt\,\Bigr|\,\F_s\Bigr\}dx.
\eaaa
\par
Let $\rho\in Z_{s}^0$ be such that \be \rho\ge 0,\quad
\int_{D}\rho(x,\o)dx=1, \quad \rho(x,\o)=0 \label{rho}\ee for all
$\o$. Let $a\in L_2(\O,\F,\P;\R^n)$ be such that $a\in D$ a.s.,  $a$
has the conditional given $\F_s$ probability density function $\rho$
on $D$, and $a$ is independent from $(w(t)-w(t_1),\w w(t)-w(t_1))$
for all $t>t_1>s$.
 Let $y(t)$ be
the solution of  Ito equation (\ref{yxs}) with initial condition
$y(s)=a$, i.e., $y(t)=y^{a,s}(t)$. In addition, let
$\g(t)=\g^{a,s}(t)$. Then
 \baaa
|(\ww p^{(k)}(\cdot,s)-\ww p(\cdot,s),\rho)_{Z_T^0}|&\le&
\E\g(T)|\Psi^{(k)}(y(T))-\Psi(y(T))|\\ &+&\E\int_s^T
\g(t)\,|\ww\xi^{(k)}(y(t),t,\o)-\xi(y(t),t,\o)|\,dt.\eaaa
\par
Let $\oo Z_{s}^0=Z_s^0$ be the space defined similarly to $Z_s^0$
but with $D$ replaced by $\R^n$. Let $u\defi \oo\L(s,T)\rho$, where
the operator $\oo\L(s,T)$ is defined similarly to $\L(s,T)$ but such
that $D$ is  replaced by $\R^n$. If $D=\R^n$, then $\oo
Z_{s}^0=Z_s^0$ and $\oo\L(s,T)=\L(s,T)$. The conditions of Theorem
5.3.1 from Rozovskii (1990)\index{, p.164} are satisfied. By this
theorem, it follows that $$ \int_{\R^n}
u(x,t,\o)\phi(x,\o)dx=\E\Big\{\g(t)\phi(y(t),\o)\,|\,\F_t\Big\}
\quad\hbox{a.s.} $$ for all $t\in [s,T]$ for any bounded function
$\phi\in\oo Z_t^0$. In fact, the cited theorem from Rozovskii (1990)
states it for non-random $\phi$, but clearly  it is also correct for
the case of $\phi\in\oo Z_t^0$ since $\phi$ is non-random
conditionally given $\F_t$.  (We can use also Theorem 2.2 from
Dokuchaev (1995)). It follows that \baaa &&|(\ww
p^{(k)}(\cdot,s)-\ww p(\cdot,s),\rho)_{Z_T^0}|\\&&\le
\E\int_{\R^n}\! u(x,T,\o)|\Psi^{(k)}(x,\o)-\Psi(x,\o)|dx
+\E\int_s^T\!
dt\int_{\R^n} u(x,t,\o)|\ww\xi^{(k)}(x,t,\o)-\xi(x,t,\o)|\,dx\\
&&\le \|u\|_{Y^1(s,T)}\left(\|\Psi^{(k)}-\Psi\|_{Z^0_T}+
\|\ww\xi^{(k)}-\xi\|_{X^{-1}}\right).\eaaa By (\ref{xixi}), it
follows that (\ref{wpp}) holds for all $\rho\in Z_{s}^0$ such that
(\ref{rho}) holds.  It follows that (\ref{wpp}) holds for any
$\rho\in Z_{s}^0$, since it can be presented as
$\rho=c_-\rho_+-c_+\rho_-$, where $\rho_\pm$ are elements of
$\Z_s^0$ such that (\ref{rho}) holds for a.e. $\o$, and $c_\pm\in\R$
are some constants.
\par This completes the proof of Theorem \ref{Threpr}
for the case when $\xi$ and $\Psi$ are bounded (and finitely
supported in $x$ if $D=\R^n$).
\par
For case of $\xi$ and $\Psi$ of the general type, it suffices to
prove theorem only when $\xi\ge 0$ and $\Psi\ge 0$. The proof for
$\xi$ and $\Psi$ with variable signs follows immediately, if we
use the linearity of (\ref{pintheorem}) and (\ref{repr})  with
respect to $(\xi,\Psi)$ and observe that $\xi=(\xi)^+-(-\xi)^+$
and $\Psi=(\Psi)^+-(-\Psi)^+$, where $(x)^+\defi \max(0,x)$.
\par
Let us consider $\xi$ and $\Psi$  such that $\xi\ge 0$ and
$\Psi\ge 0$. For $M>0$, set $$\xi_M(x,t,\o)\defi
\max(\xi(x,t,\o),M)\Ind_{\{|x|\le M\}},\quad \Psi_M(x,\o)\defi
\max(\Psi(x,\o),M)\Ind_{\{|x|\le M\}}.$$ Let $p_M\defi
L^*\xi_M+(\I_T L)^*\Psi_M$. We have proved that \baaa
p_M(x,s,\o)=\E\Big\{\g^{x,s}(T)\Psi_M(y^{x,s}(T))\Ind_{\{T\ge
\tau^{x,s}\}}\,|\,\F_s\Big\}+\E\biggl\{\int_s^{\tau^{x,s}}
\g^{x,s}(t)\,\xi_M(y^{x,s}(t),t,\o)\,dt\,\bigr|\,\F_s\biggr\}\\
\hbox{for a.e.}\ x,\o. \nonumber\eaaa By the Lebesgue's Dominated
Convergence Theorem,  it follows that  $$
\|\xi_M-\xi\|_{X^0}+\|\Psi_M-\Psi\|_{Z^0_T} \to
0\quad\hbox{as}\quad m\to +\infty. $$ By Theorem \ref{Th4.2}, it
follows that $\|p_M-p\|_{Y^1}\to 0$. On the other hand,
$\xi_M(x,t,\o)\to \xi(x,t,\o)$ and $\Psi_M(x,\o)\to\Psi(x,\o)$
from below for all $x,t,\o$ (and these sequence are non-decreasing
in $m$). Hence $p_M$ converges to the right hand part of
(\ref{repr}).
 This
completes the proof of Theorem \ref{Threpr}. $\Box$
\begin{remark} We used Theorem
  \ref{Th2FE} to obtain (\ref{X2e}) via Theorem
  \ref{Th2FE}
\end{remark}
\section{Applications: probability density for the process being
killed on the boundary} Let $s\in [0,T)$. Let $\rho\in Z_{s}^0$ be
such that $\rho\ge 0$ and $\int_{D}\rho(x,\o)dx=1$ for all $\o$. Let
$a\in L_2(\O,\F,\P;\R^n)$ be a vector such that $a\in D$ and it has
the conditional (relative to $\F_s$) probability density function
$\rho$. We assume also that $a$ is independent from $(w(t)-w(t_1),\w
w(t)-\w w(t_1))$ for all $t>t_1>s$. \par Let $u=\L(s,T)\rho$, i.e.,
$u=u(x,t,\o)$ is the solution of  the problem \be \label{4.1u} \ba
d_tu=\A u\, dt +
\sum_{i=1}^NB_iu\,dw_i(t), \quad t\ge s,\\
u|_{t=s}=\rho,\quad\quad u(x,t,\o)|_{x\in \p D}=0.\ea \ee
\par
We assume below that the assumptions of Theorem \ref{Threpr} for
$(b,\w f,\w \lambda, \b_i,\w \b_i)$ are satisfied.
\begin{theorem}\label{Thprobd} Let $s\in [0,T)$.
 Let $y(t)=y^{a,s}(t)$ be the solution of  Ito equation
(\ref{yxs}) with the initial condition $y(s)=a$. Then \be \int_{D}
u(x,T,\o)\Psi(x,\o)dx=\E\Big\{\g^{a,s}(T)\Psi(y^{a,s}(T))\Ind_{\{T\le
\tau^{a,s}\}}\,|\,\F_T\Big\} \quad\hbox{a.s.} \label{probd}\ee for
all bounded functions $\Psi\in Z_T^0$.
\end{theorem}
Note that if $D=\R^n$ then this theorem repeats Theorem 5.3.1 from
Rozovskii (1990).  However, this result is new for the case when
$D\neq \R^n$.
\begin{corollary} If $\w\b_i\equiv 0$ for all $i$ then
(\ref{probd}) means that $u(x,T,\o)$ is the conditional (relative to
$\F_T$) probability density function of the process
$y(T)=y^{a,s}(T)$ if this process is being killed at $\p D$ and if
it is being killed inside $D$ with the rate of killing $\w\lambda$.
\end{corollary}
\par
{\it Proof of Theorem \ref{Thprobd}.} It suffices to consider $s=0$
only. Let  $\Psi\in Z_T^0$ and $\w\Psi(x,\o) =\eta(\o)\w\Psi(x,\o)$,
where $\eta\in L_{\infty}(\O,\F_T,\P)$. Let $p\defi (\I_T
L)^*\w\Psi$. By Theorem \ref{Th4.2}, it follows that $$
(u(\cdot,T),\w\Psi)_{Z_T^0}=(\I_T\L\rho,
\w\Psi)_{Z_T^0}=(\rho,(\I_T\L)^*\w\Psi)_{Z_T^0}=(\rho,p(\cdot,0))_{Z_T^0}.$$
By Theorem \ref{Threpr}, \baaa (\rho,p(\cdot,0))_{Z_T^0}=
\E\int_D\rho(x)\g^{x,0}(T)\w\Psi(y^{x,0}(T)\Ind_{\{T\ge
\tau^{x,0}\}}dx= \E\eta\g^{a,0}(T)\Psi(y^{a,0}(T)\Ind_{\{T\ge
\tau^{a,0}\}}\\= \E\eta\E\{\g^{a,0}(T)\Psi(y^{a,0}(T)\Ind_{\{T\ge
\tau^{a,0}\}}|\F_T\}. \eaaa Then $$ \E\eta\int_{D}
u(x,T,\o)\Psi(x,\o)dx=
\E\eta\E\{\g^{a,0}(T)\Psi(y^{a,0}(T)\Ind_{\{T\ge
\tau^{a,0}\}}|\F_t\}.  $$ Remind that $\eta\in
L_{\infty}(\O,\F_T,\P)$ is arbitrary. Then the proof follows. $\Box$
\section{Applications: maximum principle and contraction property}
Remind that the assumptions of Theorem \ref{Threpr} for $(b,\w f,\w
\lambda, \b_i,\w \b_i)$ are satisfied.
\begin{theorem} 
\label{Th6.1} (Maximum principle) Let $\xi\in X^0$ and $\Psi\in
Z^0_T$ be such that $\xi(x,t,\o)\ge 0$ and $\Psi(x,\o)\ge 0$ for
a.e. $x,t,\o$. Then  the solution $p$ of (\ref{4.8}) is such that
$p(x,t,\o)\ge 0$ for all $t$ for a.e.  \ $t,\o$.
\end{theorem}
\par
{\it Proof.} Assume that $\xi(x,t,\o)\ge  0$ and $\Psi(x,\o)\ge 0$
for all $x,t,\o$ and that these functions have the same
measurability as described in Theorem \ref{Threpr}. In this case,
the proof follows immediately from Theorem \ref{Threpr}. Further,
let $\xi(x,t,\o)\ge 0$ and $\Psi(x,\o)\ge 0$ for a.e. $x,t,\o$.
Replace these function by some equivalent non-negative functions
$\xi'$ and $\Psi'$. Since $p=L^*\xi+(\I_T L)^*\Psi$, it follows that
$p=L^*\xi'+(\I_T L)^*\Psi'$ as an element of $Y^2$. By Theorem
\ref{Threpr} again, $p$ is nonnegative up to equivalency. Then the
proof follows. $\Box$
\begin{theorem} 
\label{Th6.2} (Maximum principle)  Let $\varphi\in X^0$ and $\Phi\in
Z_0^0$ be given such that $\varphi(x,t,\o)\ge 0$ and $\Phi(x,\o)\ge
0$ for a.e. $x,t,\o$. Then the solution $u$ of problem
 (\ref{4.1}) is such that $u(x,t,\o)\ge 0$ for all $t$ for a.e. $x,\o$.
\end{theorem}
\par
{\it Proof.} It suffices to consider $t=T$ only. Let $\Psi\in Z_T^0$
be an arbitrary function  such that $\Psi\ge 0$ a.e.  We have $$
(u(\cdot,T),\Psi)_{Z_T^0}= (\I_T
L\varphi+\I_T\L\Phi,\Psi)_{Z_T^0}=(\varphi,p)_{Z_T^0}+(\Phi,p(\cdot,0))_{Z_T^0},
$$ where $p\defi (\I_T L)^*\Psi$. Then $p(x,s,t)\ge 0$ for all $s$
for a.e. $x,\o$, and $(u(\cdot,T),\Psi)_{Z_T^0}\ge 0$. Then the
proof follows. $\Box$
\begin{theorem} 
\label{Th6.3}   (Contraction property) Under the assumptions of
Theorem \ref{Threpr}, let $\w\lambda(x,t,\o)\ge 0$ and $\w\b_i\equiv
0$ for all $i$. Then \baaa \esssup_{x,\o}|p(x,t,\o)\le
\esssup_{x,\o}|\Psi(x,\o)|+ (T-t)\esssup_{x,t,\o}|\xi(x,t,\o)|\quad
\forall t\in[0,T]. \eaaa
\end{theorem}
\par
{\it Proof.} Note that there are  bounded functions $\xi'$ and
$\Psi'$
 that are equivalent to $\xi$ and $\Psi$ and such that
 $$\esssup_{x,\o}|\Psi(x,\o)|+
(T-t)\esssup_{x,t,\o}|\xi(x,t,\o)|= \sup_{x,\o}|\Psi'(x,\o)|+
(T-t)\sup_{x,t,\o}|\xi'(x,t,\o)|. $$
  Since $p=L^*\xi+(\I_T
L)^*\Psi$, it follows that $p=L^*\xi'+(\I_T L)^*\Psi'$ as an element
of $Y^2$. It follows immediately from Theorem \ref{Threpr} that
\baaa \esssup_{x,\o}|p(x,t,\o)\le \sup_{x,\o}|\Psi'(x,\o)|+
(T-t)\sup_{x,t,\o}|\xi'(x,t,\o)|\quad \forall t\in[0,T]. \eaaa Then
the proof follows. $\Box$
\begin{theorem}
\label{Th6.4}   (Contraction property) Let $\w\lambda(x,t,\o)\ge 0$
and $\w\b_i\equiv 0$ for all $i$, and let $\varphi\in X^0$ and
$\Phi\in Z_0^0$ be given. Then the following holds for the solution
$u$ of problem (\ref{4.1}):
\begin{itemize}
\item[(a)] If $\varphi\equiv 0$, then \baaa
\E\int_D|u(x,T,\o)|dx\le \E\int_D|\Phi(x,\o)|dx. \eaaa \item[(b)] If
$\Phi=0$, then \baaa &&\E\int_D|u(x,T,\o)|dx\le
\frac{1}{T}\E\int_Q|\varphi(x,t,\o)|dxdt. \eaaa
\end{itemize}
\end{theorem}
\par
{\it Proof.} Let $\Psi\in Z_T^0$ be an arbitrary function. By
Theorem \ref{Th4.2}, it follows that  $$ (u(\cdot,T),\Psi)_{Z_T^0}=
(\I_T
L\varphi+\I_T\L\Phi,\Psi)_{Z_T^0}=(\varphi,p)_{Z_T^0}+(\Phi,p(\cdot,0))_{Z_T^0},
$$ where $p\defi (\I_T L)^*\Psi$. Then the proof follows from
Theorem \ref{Th6.3}. $\Box$

\subsection*{Conclusions}
We obtained the representation theorem  for non-Markov Ito processes
in bounded domains when the first exit times are involved. This
result is not particularly surprising;  the similar result  without
first exit times for the processes in the entire space  was obtained
long time ago. However, the setting with first exit times required
to overcome one crucial obstacle: insufficiency of the known
regularity for backward SPDEs in domains with boundaries.
Consequently, there is a little known   about  first exit times of
non-Markov processes. The representation theorem opens some further
opportunities for studying
 first exit times for non-Markov processes. It is unclear yet if it
 is possible
 to relax the strengthened coercivity required by Condition \ref{condA}.
 Probably, is some cases, this condition may be lifted  via the
 estimates from Dokuchaev (2008). To cover more general models, we suggest
to include the case of infinite number of driving Wiener processes
and more general boundary conditions. We leave it for future
research.

\subsection*{Acknowledgment}  This work  was supported by NSERC
grant of Canada 341796-2008 to the author.
\section*{References}
$\hphantom{XX}$Al\'os, E., Le\'on, J.A., Nualart, D. (1999).
 Stochastic heat equation with random coefficients
 {\it
Probability Theory and Related Fields} {\bf 115}, 1, 41-94.
\par
Bally, V., Gyongy, I., Pardoux, E. (1994). White noise driven
parabolic SPDEs with measurable drift. {\it Journal of Functional
Analysis} {\bf 120}, 484 - 510.
\par
Confortola, F. (2007).  Dissipative backward stochastic
differential equations with locally Lipschitz nonlinearity, {\it
Stochastic Processes and their Applications} {\bf 117}, Issue 5,
613-628.
\par
Chojnowska-Michalik, A., and Goldys, B. (1965). {Existence,
uniqueness and invariant measures for stochastic semilinear
equations in Hilbert spaces},  {\it Probability Theory and Related
Fields},  {\bf 102}, No. 3, 331--356.
\par
Da Prato, G., and Tubaro, L. (1996). { Fully nonlinear stochastic
partial differential equations}, {\it SIAM Journal on Mathematical
Analysis} {\bf 27}, No. 1, 40--55.
\par
Dokuchaev, N.G. (1992). { Boundary value problems for functionals
of
 Ito processes,} {\it Theory of Probability and its Applications}
 {\bf 36 }, 459-476.
\par
 Dokuchaev, N.G. (1995). Probability distributions  of  Ito's
processes: estimations for density functions and for conditional
expectations of integral functionals. {\it Theory of Probability
and its Applications} {\bf 39} (4),  662-670.
\par Dokuchaev, N.G. (2003). Nonlinear
parabolic Ito's equations and duality approach, {\it Theory of
Probability and its Applications} {\bf 48} (1), 45-62.
\par
Dokuchaev, N.G. (2005). Parabolic Ito equations and second
fundamental inequality.  {\it Stochastics} {\bf 77}, iss. 4.,
349-370.
\par
Dokuchaev, N.G. (2006). Backward parabolic Ito equations and second
fundamental inequality. Working paper: arXiv:math/0606595v3
[math.PR].
\par
Dokuchaev, N. (2008). Estimates for first exit  times of
non-Markovian It\^o processes.  {\it Stochastics} {\bf 80},
397--406.
\par Dokuchaev, N. (2010). Duality and semi-group property
for backward parabolic Ito equations. {\em Random Operators and
Stochastic Equations} {\bf 18}, 51-72.
\par
 Hu, Y., Ma, J., and Yong, J. (2002).  On Semi-linear degenerate
backward stochastic partial differential equations. {\it Probability
Theory and Related Fields}, {\bf 123}, No. 3, 381-411.
\par
Krylov, N.V.  (1980),  {\em Controlled Diffusion Processes}.
Shpringer-Verlag.
\par Gy\"ongy, I. (1998). Existence and
uniqueness results for semilinear stochastic partial differential
equations. {\it Stochastic Processes and their Applications} {\bf
73} (2), 271-299.
\par
Kim, Kyeong-Hun (2004). On stochastic partial di!erential
equations with variable coefficients in $C^1$-domains. Stochastic
Processes and their Applications {\bf 112}, 261--283.
\par Krylov, N. V. (1999). An
analytic approach to SPDEs. Stochastic partial differential
equations: six perspectives, 185--242, Math. Surveys Monogr., 64,
Amer. Math. Soc., Providence, RI.
\par
 Ladyzhenskaya, O. A., and Ural'tseva, N.N.  (1968). {\it
  Linear and quasilinear elliptic equations}. New York: Academic Press.
  \par
Ladyzhenskaia, O.A. (1985). {\it The Boundary Value Problems of
Mathematical Physics}. New York: Springer-Verlag.
\par
Ma, J., and Yong, J. (1999),  On Linear Degenerate Backward
Stochastic PDE's. {\it Probability Theory and Related Fields},
{\bf 113}, 135-170.
\par
Maslowski, B. (1995). { Stability of semilinear equations with
boundary and pointwise noise}, {\it Ann. Scuola Norm. Sup. Pisa
Cl. Sci.} (4), {\bf  22}, No. 1, 55--93.
\par
Pardoux, E. (1993).
 Stochastic partial differential equations, a review, {\it Bull. Sc. Math.}
 {\bf 117}, 29-47.
 \par
Pardoux, E., S. Peng, S. (1990). Adapted solution of a backward
stochastic differential equation. {\it System \& Control Letters}
{\bf 14}, 55-61.
\par
Pardoux, E., A. Rascanu,A. (1998). Backward stochastic
differential equations with subdifferential operators and related
variational inequalities, {\it Stochastic Process. Appl.} {\bf 76}
(2), 191-215.
\par
Rozovskii, B.L. (1990). {\it Stochastic Evolution Systems; Linear
Theory and Applications to Non-Linear Filtering.} Kluwer Academic
Publishers. Dordrecht-Boston-London.
\par
Walsh, J.B. (1986). An introduction to stochastic partial
differential equations, {\it Ecole d'Et\'e de Prob. de St.} Flour
XIV, 1984, Lect. Notes in Math 1180, Springer Verlag.
\par
Yong, J., and Zhou, X.Y. (1999). { Stochastic controls:
Hamiltonian systems and HJB equations}. New York: Springer-Verlag.
\par
Yosida, K. (1995) {\it Functional Analysis.}  Springer, Berlin
Heilderberg New York.
\par
 Zhou, X.Y. (1992). { A duality analysis on stochastic partial
differential equations}, {\it Journal of Functional Analysis} {\bf
103}, No. 2, 275--293.

\end{document}